# ALGEBRA AND LOGIC. SOME PROBLEMS

BORIS PLOTKIN

## Introduction

The paper has a form of a talk on the given topic. This second paper continues the first one [14]. It consists of three parts, ordered in a way different from that of [14]. The accents are also different. On our opinion, some simple proofs make the paper more vital.

The first part of the paper contains main notions, the second one is devoted to logical geometry, the third part describes types and isotypeness. The problems are distributed in the corresponding parts. The whole material oriented towards universal algebraic geometry (UAG), i.e., geometry in an arbitrary variety of algebras $\Theta$. We consider logical geometry (LG) as a part of UAG. This theory is strongly influenced by model theory and ideas of A.Tarski and A.I.Malcev.

I remember that A.I. Malcev, founding the journal "Algebra and logic" in Novosibirsk had in mind a natural interrelation of these topics.

Let a variety of algebras $\Theta$ be fixed. Let $W = W(X)$ be the free in $\Theta$ algebra over a set of variables $X$. The set $X$ is assumed to be finite, if the contrary is not explicitly stated. All algebras under consideration are algebras in $\Theta$. Logic is also related to the variety $\Theta$. As usual, the signature of $\Theta$ may contain constants.

## 1. Main notions

In this section we consider a system of notions, we are dealing with. Some of them are not formally defined in this paper. For the precise definitions and references use [8],[2], [9],[12],[16], [4].

The general picture of relations between these notions brings forward a lot of new problems, formulated in the following two sections. These problems are the main objective of the paper. Some results are also presented.

1.1. **Equations, points, spaces of points and algebra of formulas $\Phi(X)$.** Consider a system $T$ of equations of the form $w = w'$, $w, w' \in W(X)$.

Each system $T$ determines an algebraic set of points in the corresponding affine space over the algebra $H \in \Theta$ for every $H$ and every finite $X$.

Let $X = \{x_1, \ldots, x_n\}$. We have an affine space $H^X$ of points $\mu : X \to H$. For every $\mu$ we have also the $n$-tuple $(a_1, \ldots, a_n) = \bar{a}$ with $a_i = \mu(x_i)$. For the given $\Theta$ we have the homomorphism

$$\mu : W(X) \to H$$

and, hence, the affine space is viewed as the set of homomorphisms

$$Hom(W(X), H).$$

The classical kernel $Ker(\mu)$ corresponds to each point $\mu : W(X) \to H$.

Every point $\mu$ has also the logical kernel $LKer(\mu)$. The point is that along with the algebra $W(X)$ we will consider the algebra of formulas $\Phi(X)$. Logical kernel $LKer(\mu)$ consists of all formulas $u \in \Phi(X)$ valid on the point $\mu$.

The algebra $\Phi(X)$ will be defined later on, but let us note now that it is an extended Boolean algebra (Boolean algebra in which quantifiers $\exists x, x \in X$ act as operations, and equalities ($\Theta$-equalities) $w \equiv w'$, $w, w' \in W(X)$ are defined). It is also defined what does it mean that the point $\mu$ satisfies a formula $u \in \Phi(X)$. These $u$ are treated as equations. For $T \subset \Phi(X)$ in $Hom(W(X), H)$ we have an elementary set (definable set) consisting of points $\mu$ which satisfy every $u \in T$.

Each kernel $LKer(\mu)$ is a Boolean ultrafilter in $\Phi(X)$. Note that

$$Ker(\mu) = LKer(\mu) \cap M_X,$$

where $M_X$ is the set of all $w \equiv w'$, $w, w' \in W(X)$.

**1.2. Extended Boolean algebras.** Let us make some comments regarding the definition of the notion of extended Boolean algebra.

Let $B$ be a Boolean algebra. An existential quantifier on $B$ is an unary operation $\exists : B \to B$ subject to conditions

(1) $\exists(0) = 0$,
(2) $a \leq \exists(a)$,
(3) $\exists(a \wedge \exists b) = \exists a \wedge \exists b$.

The *universal* quantifier $\forall : B \to B$ is defined dually:

(1) $\forall(1) = 1$,
(2) $a \geq \forall(a)$,
(3) $\forall(a \vee \forall b) = \forall a \vee \forall b$.

Here the numerals 0 and 1 are zero and unit of the Boolean algebra $B$ and $a, b$ are arbitrary elements of $B$.

As usual, the quantifiers $\exists$ and $\forall$ are coordinated by: $\neg(\exists a) = \forall(\neg a)$, and $(\forall a) = \neg(\exists(\neg a))$.

Now suppose that a variety of algebras $\Theta$ is fixed and $W(X)$ is the free in $\Theta$ algebra over the set of variables $X$. These data allow to define the extended Boolean algebra. This is a Boolean algebra where the quantifiers $\exists x$ are defined for every $x \in X$ and

$$\exists x \exists y = \exists y \exists x,$$

for every $x$ and $y$ from X. Besides that, for every pair of elements $w, w' \in W(X)$ in an extended Boolean algebra the equality $w \equiv w'$ is defined. These equalities are considered as nullary operations, that is as constants. Each equality satisfies conditions of an equivalence relation, and for every operation $\omega$ from the signature of algebras from $\Theta$ we have

$$w_1 \equiv w'_1 \wedge \ldots \wedge w_n \equiv w'_n \to w_1 \ldots w_n \omega \equiv w'_1 \ldots w'_n \omega.$$

Note that quantifiers and Boolean connectives satisfy

$$\exists(a \vee b) = \exists a \vee \exists b$$

$$\forall(a \wedge b) = \forall a \wedge \forall b$$

Algebra of formulas $\Phi(X)$ is an example of extended Boolean algebra in $\Theta$. Now consider another example.

1.3. **Important example.** Let us start from an affine space $Hom(W(X), H)$. Let $Bool(W(X), H)$ be a Boolean algebra of all subsets in $Hom(W(X), H)$. Extend this algebra by adding quantifiers $\exists x$ and equalities. For $A \in Bool(W(X), H)$ we set: $B = \exists x A$ is a set of points $\mu : W(X) \to H$ such that there is $\nu : W(X) \to H$ in $A$ and $\mu(x') = \nu(x')$ for $x' \in X$, $x' \neq x$. It is indeed an existential quantifier for every $x \in X$.

Define an equality in $Bool(W(X), H)$ for $w \equiv w'$ in $M_X$. Denote it by $[w \equiv w']_H$ and define it, setting $\mu \in [w \equiv w']_H$ if $(w, w') \in Ker(\mu)$, i.e., $w^\mu = w'^\mu$.

**Remark 1.1.** The set $[w \equiv w']_H$ can be empty. Thus we give the following definition. The equality $[w \equiv w']_H$ is called admissible for the given $\Theta$ if for every $H \in \Theta$ the set $[w \equiv w']_H$ is not empty. If $\Theta$ is the variety of all groups, then each equality is admissible. The same is true for the variety of associative algebras with unity over complex numbers. However, for the field of real numbers this is not the case. Here $x^2 + 1 = 0$ is not an admissible equality.

We assume that in each algebra of formulas $\Phi(X)$ lie all $\Theta$-equalities. To arbitrary equality $w \equiv w'$ corresponds either a non-empty equality $[w \equiv w']_H$ in $H \in \Theta$, or the empty set in $H \in \Theta$ which is the zero element of this boolean algebra.

We arrived to an extended Boolean algebra, denoted now by $Hal_\Theta^X(H)$. This algebra and the algebra of formulas $\Phi(X)$ have the same signature.

1.4. **Homomorphism** $Val_H^X$. We will proceed from the homomorphism
$$Val_H^X : \Phi(X) \to Hal_\Theta^X(H),$$
with the condition $Val_H^X(w \equiv w') = [w \equiv w']_H$ for equalities, if $[w \equiv w']_H$ non-empty, or 0 otherwise. This homomorphism will be defined in subsection 1.9. The existence of such homomorphism is not a trivial fact, since the equalities $M_X$ does not generate (and of course does not generate freely) the algebra $\Phi(X)$. If, further, $u \in \Phi(X)$, then $Val_H^X(u)$ is a set of points in the affine space $Hom(W(X), H)$. We say that a point $\mu$ satisfies the formula $u$ if $\mu$ belongs to $Val_H^X(u)$. Thus, $Val_H^X(u)$ is precisely the set of points satisfying the formula $u$. Define the logical kernel $LKer(\mu)$ of a point $\mu$ as a set of all formulas $u$ such that $\mu \in Val_H^X(u)$.

We have also
$$Ker(\mu) = LKer(\mu) \cap M_X.$$

Here $Ker(\mu)$ is the set of all formulas of the form $w = w'$, $w, w' \in W(X)$, such that the point $\mu$ satisfies these formulas. In parallel, $LKer(\mu)$ is the set of all formulas $u$, such that the point $\mu$ satisfies these formulas.

Note that that in Section 3.7 we consider the notion of the type $Tp^H(\mu)$. It is also a sort of a kernel of the point $\mu$. It consists of $X$-special formulas $u$ which are satisfied by $\mu$.

Then,
$$Ker(Val_H^X) = Th^X(H),$$
$$\bigcap_{\mu:W(X) \to H} LKer(\mu) = Th^X(H).$$

Here $Th^X(H)$ is a set of formulas $u \in \Phi(X)$, such that $Val_H^X(u)$ is unit in $Bool(W(X), H)$. That is $Val_H^X(u) = Hom(W(X), H)$ and thus $Th^X(H)$ is an $X$-component of the elementary theory of the algebra $H$.

In general we have a multi-sorted representation of the elementary theory
$$Th(H) = (Th^X(H), X \in \Gamma^0).$$

It follows from the previous considerations that the algebra of formulas $\Phi(X)$ can be embedded in $Hal_\Theta^X(H)$ modulo elementary theory of the algebra $H$. This fact will be used in the sequel.

### 1.5. Multi-sorted logic: first approximation.

Let, further, $X^0$ be an infinite set of variables and $\Gamma^0$ a system of all finite subsets $X$ in $X^0$.

So, in the logic under consideration we have an infinite system $\Gamma^0$ of finite sets instead of one infinite $X^0$. This leads to multi-sorted logic. This approach is caused by relations with UAG. We distinguish UAG, equational UAG and LG in UAG. Correspondingly, we have algebraic sets of points and elementary sets of points in the affine space. In the third part of the paper along with the system of sorts $\Gamma^0$ we use also a system of sorts $\Gamma$ where one initial infinite set $X^0$ is added to the system $\Gamma^0$.

### 1.6. Algebra $Hal_\Theta(H)$.

Define now an important category $Hal_\Theta(H)$. Denote first a category of all $Hom(W(X), H)$, $X \in \Gamma$ by $\Theta(H)$. For every $s : W(X) \to W(Y)$ consider a mapping

$$\tilde{s} : Hom(W(Y), H) \to Hom(W(X), H)$$

defined by the rule $\tilde{s}(\mu) = \mu s : W(X) \to H$ for $\mu : W(Y) \to H$. These $\tilde{s}$ are morphisms of the category $\Theta(H)$. We have a contra-variant functor $\Theta^0 \to \Theta(H)$. It is proved [12] that this functor determines duality of the categories if and only if $\Theta = Var(H)$.

Define further a category of all $Hal_\Theta^X(H)$. Proceed from $s : W(X) \to W(Y)$ and pass to $\tilde{s} : Hom(W(Y), H) \to Hom(W(X), H)$. Take

$$s_*(A) = \tilde{s}^{-1}(A) = B \subset Hom(W(Y), H)$$

for every object $A \subset Hom(W(X), H)$. We have $\mu \in B$ if and only if $\mu s = \tilde{s}(\mu) \in A$. This determines a morphism

$$s_* = s_*^H : Hal_\Theta^X(H) \to Hal_\Theta^Y(H).$$

Here $s_*$ is well coordinated with the Boolean structure, and relations with quantifiers and equalities are coordinated by (definite) identities [16], [14]. The category $Hal_\Theta(H)$ can be also treated as a multi-sorted algebra

$$Hal_\Theta(H) = (Hal_\Theta^X(H), X \in \Gamma).$$

### 1.7. Variety of Halmos algebras $Hal_\Theta$.

Algebras in $Hal_\Theta$ have the form

$$\mathfrak{G} = (\mathfrak{G}_X, \ X \in \Gamma^0).$$

Here all domains $\mathfrak{G}_X$ are $X$-extended Boolean algebras. The unary operation
$$s_* : \mathfrak{G}_X \to \mathfrak{G}_Y$$
corresponds to each homomorphism $s : W(X) \to W(Y)$. Besides, we will define a category $\mathfrak{G}$ of all $\mathfrak{G}_X$, $X \in \Gamma^0$ with morphisms $s_* : \mathfrak{G}_X \to \mathfrak{G}_Y$. The transition $s \to s_*$ determines a covariant functor $\Theta^0 \to \mathfrak{G}$. Every morphism $s_*$ respects Boolean operations in $\mathfrak{G}_X$ and $\mathfrak{G}_Y$. Correlations of $s_*$ with equalities and quantifiers are described by identities (see [16], [14]). Operations of $s_*$ type make logics dynamical.

Each algebra $Hal_\Theta(H)$ belongs the variety $Hal_\Theta$. Moreover, all $Hal_\Theta(H)$, where $H$ runs $\Theta$, generate the variety $Hal_\Theta$, ( see [13]).

Recall, that every ideal of an extended Boolean algebra is a Boolean ideal invariant with respect to the universal quantifiers action. Extended Boolean algebra is called simple if it does not have non-trivial ideals. In the many-sorted case an ideal is a system of one-sorted ideals which respects all morphism of type $s_*$. A many-sorted Halmos algebra is simple if it does not have non-trivial ideals. Algebras $Hal_\Theta(H)$ and their subalgebras are simple Halmos algebras, see [17]. Moreover, these algebras are the only simple algebras in the variety $Hal_\Theta$. Finally, every Halmos algebra is residually simple, see [17]. This fact is essential in the next subsection. Note, that all these facts are true because of the clever choice of the identities in the variety $Hal_\Theta$.

1.8. **Multi-sorted algebra of formulas.** We shall define the algebra of formulas
$$\widetilde{\Phi} = (\Phi(X),\ X \in \Gamma^0).$$

We define this algebra as the free over the multi-sorted set of equalities
$$M = (M_X, X \in \Gamma^0)$$
algebra in $Hal_\Theta$. Assuming this property denote it as
$$Hal_\Theta^0 = (Hal_\Theta^X,\ X \in \Gamma^0).$$
So, $Hal_\Theta^X = \Phi(X)$ and $\widetilde{\Phi} = Hal_\Theta^0$.

In order to define $Hal_\Theta^0$ we start from the absolutely free over the same $M$ algebra
$$\mathfrak{L}^0 = (\mathfrak{L}^0(X),\ X \in \Gamma^0).$$
This free algebra is considered in the signature of the variety $Hal_\Theta$. Algebra $\mathfrak{L}^0$ can be viewed as the algebra of pure formulas of the corresponding logical calculus.

Then, $\widetilde{\Phi}$ is defined as the quotient algebra of $\mathfrak{L}^0$ modulo the verbal congruence of identities of the variety $Hal_\Theta$. The same algebra $\widetilde{\Phi}$ can be obtained from $\mathfrak{L}^0$ using the Lindenbaum-Tarski approach. Namely, basing on identities of $Hal_\Theta$ we distinguish in $\mathfrak{L}^0$ a system of axioms and rules of inference. For every $X \in \Gamma^0$ consider the formulas

$$(u \to v) \wedge (v \to u),$$

where $u, v \in \mathfrak{L}^0(X)$. Here $u \to v$ means $\neg u \vee v$. We assume that every

$$(u \to v) \wedge (v \to u),$$

is deducible from the axioms if and if the pair $(u, v)$ belongs to the $X$-component of the given verbal congruence.

So, $\widetilde{\Phi}$ can be viewed as an algebra of the compressed formulas modulo this congruence.

1.9. **Homomorphism $Val_H$.** Proceed from the mapping

$$M_X \to Hal_\Theta^X(H),$$

which takes the equalities $w \equiv w'$ in $M_X$ to the corresponding equalities $[w \equiv w']_H$ in $Hal_\Theta^X(H)$. This gives rise also to the multi-sorted mapping

$$M = (M_X, X \in \Gamma^0) \to Hal_\Theta(H) = (Hal_\Theta^X(H), \ X \in \Gamma^0).$$

Since the multi-sorted set $M$ generates freely the algebra $\widetilde{\Phi}$ this mapping is uniquely extended up to the homomorphism

$$Val_H : \widetilde{\Phi} \to Hal_\Theta(H).$$

Note that this homomorphism is a unique homomorphism from $\widetilde{\Phi} \to Hal_\Theta(H)$, since equalities are considered as constants.

We have

$$Val_H^X : \Phi(X) \to Hal_\Theta^X(H),$$

i.e., $Val_H$ acts componentwise for each $X \in \Gamma^0$.

Recall that for every $u \in \Phi(X)$ the corresponding set $Val_H^X(u)$ is a set of points $\mu : W(X) \to H$ satisfying the formula $u$ (see Subsection 1.4). The logical kernel $LKer(\mu)$ was defined in Subsection 1.1 in these terms. Now we can say, that if a formula $u$ belongs to $\Phi(X)$ and a point $\mu : W(X) \to H$ is given, then

$$u \in LKer(\mu) \text{ if and only if } \mu \in Val_H^X(u).$$

We shall note that a formula $u$ can be, in general, of the form $u = s_*(v)$, where $v \in \Phi(Y)$, $Y$ is different from $X$. This means that the logical kernel of the point is very big and it gives a rich characterization of the whole theory.

Recall further that $LKer(\mu)$ is a Boolean ultrafilter containing the elementary theory $Th^X(H)$. Any ultrafilter with this property will be considered as an $X$-type of the algebra $H$.

It is clear that
$$Ker(Val_H) = Th(H).$$
This remark is used, for example, in

**Definition 1.2. An algebra $H \in \Theta$ is called saturated if for every $X \in \Gamma$ for each ultrafilter $T$ in $\Phi(X)$ containing $Th^X(h)$ there is a representation $T = LKer(\mu)$ for some $\mu : W(X) \to H$.**

Recall that the algebra $\widetilde{\Phi}$ is residually simple. This fact implies two important observations:

1. Let $u, v$ be two formulas in $\Phi(X)$. These formulas coincide if and only if for every algebra $H \in \Theta$ the equality
$$Val_H^X(u) = Val_H^X(v)$$
holds.

2. Let a morphism $s : W(X) \to W(Y)$ be given. It corresponds the morphism $s_* : \Phi(X) \to \Phi(Y)$. Let us take formulas $u \in \Phi(X)$ and $v \in \Phi(Y)$. The equality
$$s_*(u) = v$$
holds true if and only if for every algebra $H$ in $\Theta$ we have
$$s_*(Val_H^X(u) = Val_H^Y(v).$$

The following commutative diagram relates syntax with semantics

$$\begin{array}{ccc} \Phi(X) & \xrightarrow{s_*} & \Phi(Y) \\ {\scriptstyle Val_H^X}\downarrow & & \downarrow{\scriptstyle Val_H^Y} \\ Hal_\Theta^X(H) & \xrightarrow{s_*^H} & Hal_\Theta^Y(H) \end{array}$$

We finished the survey of the notions of multi-sorted logic needed for UAG and in the next section we will relate these notions with the ideas of one-sorted logic used in Model Theory. Note also that we cannot define algebras of formulas $\Phi(X)$ individually. They are defined only in the multi-sorted case of algebras $\widetilde{\Phi} = (\Phi(X), X \in \Gamma^0)$.

In fact, the definition of the algebra of formulas $\widetilde{\Phi}$ and the system of algebras $\Phi(X)$ is the main result of the first part of the paper. They are essentially used throughout the paper.

## 2. Logical geometry

2.1. **Introduction.** The theory under consideration is universal. Its system of notions assumes arbitrary variety of algebras $\Theta$ and arbitrary algebra $H \in \Theta$. This system leads to various new problems, which are rather general as well. Concretization of varieties and algebras starts from solutions of problems. This general perspective is interesting even for the classical AG which is related to the variety $\Theta = Com - P$, the variety of all commutative and associative algebras with the unit over the field $P$.

In this section we define logically homogeneous algebras and logically perfect varieties $\Theta$. We do not know whether $\Theta = Com - P$ is logically perfect and whether a finitely generated free in $Com - P$ algebra of polynomials is logically homogeneous.

Let, further, $H_1$ and $H_2$ be two algebras in $\Theta$. The notions of geometrical equivalence and logical equivalence are defined for algebras in $\Theta$. If we have geometrical equivalence, then the corresponding categories of algebraic sets $K_\Theta(H_1)$ and $K_\Theta(H_2)$ are isomorphic. Logical equivalence implies isomorphism of categories of elementary sets $LK_\Theta(H_1)$ and $LK_\Theta(H_2)$. We have sufficient conditions of isomorphism of categories for any $\Theta$. However, if we are interested in necessary and sufficient conditions, then the results are different for different $\Theta$. Such conditions are well described for AG in $Com - P$, see [12].

Namely, let $\sigma$ be an automorphism of the field $P$, $H_1^\sigma$ - the corresponding $\sigma$-twisted algebra. The categories of algebraic sets $K_\Theta(H_1)$ and $K_\Theta(H_2)$ are isomorphic if and only if there is an automorphism $\sigma \in Aut(P)$, such that the algebras $H_1^\sigma$ and $H_2$ are AG-equivalent (see Definition 2.1) For LG the similar question is open.

Let us note the following problem for $\Theta = Com - P$:

**Problem 1. Given algebras $H_1$ and $H_2$, find necessary and sufficient conditions when the corresponding categories of elementary sets are isomorphic.**

We speak here of elementary sets in the LG case and of definable sets in the MT case (see Subsection 3.5 for definitions).

There are also other problems for LG in the variety $Com - P$. Therefore let us mention a general one:

**Problem 2.** Develop general methods for solutions of $LG$-problems in $Com - P$.

This section is devoted to the logical geometry in various $\Theta$ and, in particular, in $\Theta = Com - P$.

2.2. **Galois correspondence in the Logical Geometry.** Now, let us pass to general definitions.

Consider first the general Galois correspondence. Let $T$ be a system of equations of the $w = w'$ type, $w, w' \in W(X)$, $X \in \Gamma^0$. We set

$$A = T'_H = \{\mu : W(X) \to H \mid T \subset Ker(\mu)\}.$$

Here $A$ is an algebraic set in $Hom(W(X), H)$, determined by the set $T$.

Let, further, $A$ be a subset in $Hom(W(X), H)$. We set
$$T = A'_H = \bigcap_{\mu \in A} Ker(\mu).$$

Here $T$ is an $H$-closed congruence in $W(X)$. We have also closures $T''_H$ and $A''_H$.

Let us pass to logical geometry. Let $T$ be a set of formulas in $\Phi(X)$. We set
$$A = T^L_H = \{\mu : W(X) \to H \mid T \subset LKer(\mu)\}.$$

We have also
$$A = \bigcap_{u \in T} Val^X_H(u).$$

Here $A$ is an elementary set in $Hom(W(X), H)$, determined by the set $T$. For the set of points $A$ we set
$$T = A^L_H = \bigcap_{\mu \in A} LKer(\mu).$$

We have also $u \in T$ if and only if $A \subset Val^X_H(u)$. Here $T$ is a Boolean filter in $\Phi(X)$ determined by the set of points $A$.

2.3. **Category of elementary sets for a given algebra $H$.** Define now a category of algebraic sets $K_\Theta(H)$ and a category of elementary sets $LK_\Theta(H)$.

Define first a category $Set_\Theta(H)$. Its objects are pairs $(X, A)$ with $A$ a subset in $Hom(W(X), H)$ and $X \in \Gamma^0$. We start the definition of morphisms
$$[s] : (X, A) \to (Y, B)$$
with $s : W(Y) \to W(X)$. Further, take $\tilde{s} : Hom(W(X), H) \to Hom(W(Y), H)$. Take now those $s$ for which $\tilde{s}(\mu) = \nu \in B$ if $\mu \in A$. This defines $[s]$ and the transition $s \to [s]$ is a contravariant functor $\Theta^0 \to Set_\Theta(H)$.

Now, $K_\Theta(H)$ is a full subcategory in $Set_\Theta(H)$, whose objects $A$ are algebraic sets.

If for $A$ we take elementary sets, then we have the category $LK_\Theta(H)$ which is a full subcategory in $Set_\Theta(H)$.

Let us formulate two key definitions and two results (see, for example, [12], [15]).

**Definition 2.1. Algebras $H_1$ and $H_2$ are AG-equivalent, if $T''_{H_1} = T''_{H_2}$ always holds true.**

**Definition 2.2.** Algebras $H_1$ and $H_2$ are **LG-equivalent**, if $T_{H_1}^{LL} = T_{H_2}^{LL}$ always holds true.

**Theorem 2.3.** If $H_1$ and $H_2$ are AG-equivalent, then the categories $K_\Theta(H_1)$ and $K_\Theta(H_2)$ are isomorphic.

**Theorem 2.4.** If $H_1$ and $H_2$ are LG-equivalent, then the categories $LK_\Theta(H_1)$ and $LK_\Theta(H_2)$ are isomorphic.

**Remark 2.5.** The geometric notion of $LG$-equivalence of algebras is equivalent to modal theoretic notion of isotypic algebras (Section 3). Thus, if algebras $H_1$ and $H_2$ are isotypic, then the categories of definable sets $LK_\Theta(H_1)$ and $LK_\Theta(H_2)$ are isomorphic for every $\Theta$. On the other hand, for logically perfect varieties $\Theta$ (see Subsection 2.4 for the definition) if $H_1$ is a free in $\Theta$ algebra of the finite rank $H_2$ is any other algebra then their isotypeness already implies isomorphism of these algebras themselves.

2.4. **Logically perfect and logically regular varieties.** Let $H$ be an algebra in $\Theta$.

**Definition 2.6.** Algebra $H$ is called logically homogeneous if for every two points $\mu : W(X) \to H$ and $\nu : W(X) \to H$ the equality $LKer(\mu) = LKer(\nu)$ holds if and only if there exists an automorphism $\sigma$ of the algebra $H$ such that $\mu = \nu\sigma$.

**Definition 2.7.** A variety of algebras $\Theta$ is called logically perfect if every finitely generated free in $\Theta$ algebra $W(X)$, $X \in \Gamma^0$ is logically homogeneous.

It follows from Definition 2.2 that two algebras $H_1$ and $H_2$ are $LG$-equivalent if and only if for every $\mu : W(X) \to H_1$ there exists $\nu : W(X) \to H_2$ such that $LKer\mu = LKer\nu$.

**Definition 2.8.** An algebra $H$ in $\Theta$ is called logically separable, if every $H' \in \Theta$ which is $LG$-equivalent to $H$ is isomorphic to $H$.

**Definition 2.9.** A variety $\Theta$ is called logically regular if every free in $\Theta$ algebra $W(X)$, $X \in \Gamma^0$ is logically separable.

We will see (Theorem 3.11) that the following fact takes place:

If the variety $\Theta$ is logically perfect then $\Theta$ is logically regular.

**Problem 3.** Is the converse statement true?

**Problem 4.** Is the variety Com-P logically perfect?

2.5. **Conditions of isomorphism of categories of elementary sets.** We are interested in the following problem

**Problem 5.** Find necessary and sufficient conditions on algebras $H_1$ and $H_2$ in $\Theta$ that provide isomorphism of the categories $LK_\Theta(H_1)$ and $LK_\Theta(H_2)$.

The similar problem had been solved for the categories of algebraic sets [10], [11]. Namely,

**Definition 2.10.** Two algebras $H_1$ and $H_2$ in $\Theta$ are called geometrically similar, if there exists a commutative diagram

$$\begin{array}{ccc} \Theta^0 & \xrightarrow{\varphi} & \Theta^0 \\ & \searrow Cl_{H_1} \quad Cl_{H_2} \swarrow & \\ & PoSet & \end{array}$$

where $\varphi$ is an automorphism of $\Theta^0$, $Cl_H(W)$ is the poset of all $H$-closed congruences $T$ in $W$. Commutativity of this diagram means that there exists an isomorphism of functors $Cl_{H_1}$ and $Cl_{H_2} \cdot \varphi$. Denote this isomorphism by $\alpha(\varphi)$.

**Theorem 2.11.** *Let $Var(H_1) = Var(H_2) = \Theta$. Categories $K_\Theta(H_1)$ and $K_\Theta(H_2)$ are isomorphic if and only if the algebras $H_1$ and $H_2$ are geometrically similar.*

In view of this theorem we proceed also in the situation of logical geometry.

**Definition 2.12.** We call algebras $H_1$ and $H_2$ in $\Theta$ logically similar if there exists a commutative diagram

$$\begin{array}{ccc} \tilde{\Phi} & \xrightarrow{\varphi} & \tilde{\Phi} \\ & \searrow Cl_{H_1} \quad Cl_{H_2} \swarrow & \\ & Lat & \end{array}$$

Here $\tilde{\Phi} = (\Phi(X), X \in \Gamma^0) = Hal_\Theta^0$, and $\varphi : \tilde{\Phi} \to \tilde{\Phi}$ is an automorphism of categories. Denote the lattice of $H$-closed filters in $\Phi(X)$ by $Cl_H(\Phi(X))$ for each $X \in \Gamma^0$. Here $Cl_{H_1}$ and $Cl_{H_2}$ are corresponding functors.

Commutativity of the diagram means that there is an isomorphism of functors
$$\alpha(\varphi) : Cl_{H_1} \to Cl_{H_2}\varphi.$$

Denote $\varphi(\Phi(X)) = \Phi(Y)$. If $T$ is an $H_1$-closed filter in $\Phi(X)$, then $\alpha(\varphi)(T) = T^*$ is an $H_2$ closed filter in $\Phi(Y)$. Let now $(X, A)$ be an object of the category $LK_\Theta(H_1)$. We assume that $A = T^L_{H_1}$, where $T$ is an $H_1$ closed filter in $\Phi(X)$.

Let $(Y, B)$ be the corresponding object in the category $LK_\Theta(H_2)$. Then $B = T^{*L}_{H_2}$.

So we know how the isomorphism between $LK_\Theta(H_1)$ and $LK_\Theta(H_2)$ we are looking for, acts on objects. However, the definition of logical similarity does not allow to determine action of the needed isomorphism on morphisms. Some additional information is needed. For example, we need an information about automorphisms $\varphi$ of $\tilde{\Phi}$ and, possibly, we need to demand that $Var(H_1) = Var(H_2) = \Theta$.

So, we formulate

**Problem 6.** Find additional conditions on algebras $H_1$ and $H_2$, such that the categories $LK_\Theta(H_1)$ and $LK_\Theta(H_2)$ are isomorphic if and if $H_1$ and $H_2$ are logically similar.

2.6. **Logically noetherian and saturated algebras.**

**Definition 2.13. An algebra $H$ is called logically noetherian if for any set of formulas $T \subset \Phi(X)$, $X \in \Gamma^0$ there is a finite subset $T_0$ in $T$ determining the same set of points $A$ that is determined by the set $T$.**

**Definition 2.14. An algebra $H \in \Theta$ is called $LG$-saturated if for every $x \in \Gamma$ for each ultrafilter $T$ in $\Phi(x)$ containing $Th^x(h)$ there is a representation $T = LKer(\mu)$ for some $u : W(X) \to H$.**

**Theorem 2.15. If an algebra $H$ is logically noetherian then $H$ is LG-saturated.**

*Proof.* We start from the homomorphism:
$$Val^X_H : \Phi(X) \to Hal^X_\Theta(H).$$
Here $Ker(Val^X_H) = Th^X(H)$. Consider the quotient algebra $\Phi(X)/Th^X(H)$ which is isomorphic toi a subalgebra in $Hal^X_\Theta(H)$. For every $u \in \Phi(X)$ denote by $[u]$ the image of $u$ in the quotient algebra. By definition $[u] = 0$ means that $Val^X_H(u)$ is the empty subset in $Hom(W(X), H)$. Analogously $[u] = 1$ means that $Val^X_H(u)$ is the whole space $Hom(W(X), H)$ and, thus, $u \in Th^X(H)$.

Denote by $T$ an ultrafilter in $\Phi(X)$, containing the theory $Th^X(H)$. We need to check that there is a point $\mu : W(X) \to H$ such that $T = LKer(\mu)$. Let $[u] = 0$. Then $[\neg u] = 1$, which means that $\neg u \in Th^X(H) \subset T$. Hence $\neg \in T$. Then $u$ does not belong to $Th^X(H)$, since $T$ cannot contain both $u$ and $\neg u$. So $u \notin T$. Thus, if $[u] = 0$ then $u \notin T$. If $u \in T$, then $[u] \neq 0$. This means that $Val^X_H(u)$ is not empty. Thus, we have a point $\mu : W(X) \to H$ which satisfies $u$, that is $u \in LKer(\mu)$. Since $H$ is logically noetherian then there exists a finite subset $T_0 = \{u_1, \ldots, u_n\}$ such that $T_H^L = (T_0)_H^L$. Take $u = u_1 \wedge u_2 \wedge \ldots u_n$. Since all $u_i \in T$ then $u \in T$, Then there exists $\mu$ satisfying formula $u$. The same point $\mu$ satisfies every $u_i$. Hence $\mu \in (T_0)^L_( H) = T^L_( H)$ and $T$ lies in $LKer(\mu)$. Therefore $T = LKer(\mu)$. $\square$

Each finite algebra $H$ is logically noetherian. Therefore, every finite $H$ is saturated. This holds for every $\Theta$.

**Problem 7. Find other interesting saturated algebras $H$ in various varieties $\Theta$.**

We proceed from logical geometry approach to the saturation of algebra problem. It is clear that there is also model theoretic approach.

**Problem 8. Whether it is true that both approaches lead to one and the same saturation of an algebra.**

**Problem 9. Treat separately the idea of saturated algebra in the variety of abelian groups.**

We have already mentioned that the group $Aut(H)$ acts in each space $Hom(W(X), H)$, $X \in \Gamma$.

**Definition 2.16.** Let us call an algebra $H$ **automorphically finitary** if in each such action there is only finite number of $Aut(H)$-orbits.

It is easy to show that if algebra $H$ is automorphically finitary, then it is logically noetherian. There exist infinite automorphically finitary algebras and, thus, there are infinite saturated algebras.

There are also saturated abelian groups.

**Problem 10. Describe all automorphically finitary abelian groups.**

**Problem 11. Consider examples of non-commutative automorphically finitary groups.**

**Problem 12. Classify abelian groups by $LG$-equivalence relation.**

We used here some notions which will be discussed in more detail in the next section.

We had considered two important characteristics of varieties of algebras, namely, their logical perfectness and logical regularity. Let us introduce one more characteristic.

We call a variety $\Theta$ exceptional if

- any two free in $\Theta$ algebras $W(X)$ and $W(Y)$ of a finite rank, generating the whole $\Theta$, are elementary equivalent, and
- if they are isotypic then they are isomorphic.

**Problem 13. Whether it is true that only the variety of all groups is unique.**

## 3. Types and isotypeness

3.1. **Definitions of types.** The notion of a type is an important notion of Model Theory. We will distinguish model theoretical types (MT-types) and Logical Geometry types (LG-types). Both kinds of types are oriented towards some algebra $H \in \Theta$. Model-theoretical idea of a type is described in many sources, see, in particular, [3]. A type of a point $\mu : W(X) \to H$ is a logical characteristic of the point. We consider this idea from the perspective of algebraic logic (cf., [16]).

All the definitions are given in the terms of algebraic logic. Proceed from the algebra of formulas $\Phi(X^0)$. It arrives from the algebra of pure first-order formulas with equalities $w \equiv w'$, $w, w' \in W(X^0)$ by Lindenbaum-Tarski algebraization principle. $\Phi(X^0)$ is an $X^0$-extended Boolean algebra, which means that $\Phi(X^0)$ is a boolean algebra with quantifiers $\exists x$, $x \in X^0$ and equalities $w \equiv w'$, where $w, w' \in W(X^0)$. All these equalities generate the algebra $\Phi(X^0)$. Besides, the semigroup $End(W(X^0))$ acts in the boolean algebra $\Phi(X^0)$ and we can speak of a polyadic algebra $\Phi(X^0)$. However, the elements $s \in End(W(X^0))$ and the corresponding $s_*$ are not included in the signature of the algebra $\Phi(X^0)$.

For the formulas $u \in \Phi(X^0)$ one can speak, as usual, about free and bound occurrences of the variables from $u$.

Define further $X$-special formulas in $\Phi(X^0)$, $X = \{x_1, \ldots, x_n\}$. Take $X^0 \setminus X = Y^0$. A formula $u \in \Phi(X^0)$ is $X$-special if each its free variable is occurred in $X$ and each bound variable belongs to $Y^0$. A formula $u \in \Phi(X^0)$ is closed if it does not have free variables.

Note that denoting an $X$-special formula $u$ as $u = u(x_1, \ldots, x_n; y_1, \ldots, y_m)$ we mean solely that the set $X$ consists of variables $x_i$, $i = 1, \ldots n$, and those of them who occur in $u$, occur freely.

$X$-type is a set of $X$-special formulas, consistent with the elementary theory of the algebra $H$. Here a type is related to the algebra $H$. It is an $X$-MT-type (Model Theoretic type). As for an $X$-LG-type (Logical Geometric type), it is a boolean ultrafilter in the corresponding $\Phi(X)$, which contains the elementary theory $Th^X(H)$. An MT-type of a point $\mu : W(X) \to H$ we denote by $Tp^H(\mu)$ while LG-type of the same point is $LKer(\mu)$.

**Definition 3.1.** Let a point $\mu : W(X) \to H$, with $a_i = \mu(x_i)$, be given. An $X$-special formula $u = u(x_1, \ldots, x_n; y_1, \ldots, y_m)$ belongs to the type $Tp^H(\mu)$ if the formula $u(a_1, \ldots, a_n; y_1, \ldots, y_m)$ is satisfied in the algebra $H$.

In view of this definition a formula $u = u(x_1, \ldots, x_n; y_1, \ldots, y_m)$ is satisfied on the point $\mu$ if $u(a_1, \ldots, a_n; y_1, \ldots, y_m)$ is satisfied in the algebra $H$. Hence, $Tp^H(\mu)$ consists of all $X$-special formulas satisfied on $\mu$.

The formula $u(a_1, \ldots, a_n; y_1, \ldots, y_m)$ is closed. Thus if it is satisfied one a point then its value set $Val_X^H(u)$ is the whole $Hom(W(X), H)$.

Note that our definition of an $MT$-type differs from the standard one. In the standard definition all variables from $X = \{x_1, \ldots, x_n\}$ are exactly the set of all free variables in $u = u(x_1, \ldots, x_n; y_1, \ldots, y_n)$. In our definition the set of free variables in the formula $u$ can be a part of $X$. In particular, the set of free variables can be empty. In this case the formula $u$ belongs to the type if it is satisfied in $H$.

This insight on types is needed in order to relate $MT$-types and $LG$-types.

One more remark. When we say that a variable occur in a formula $u \in \Phi(X^0)$, this means that it occur in one of the equalities $w = w'$, participating in $u$. The set of variables occurring in $u$ determines a subalgebra $\Phi(X \cup Y)$ in $\Phi(X^0)$, such that $u \in \Phi(X \cup Y)$. If we stay in one-sorted logic, this is a subalgebra in the signature of the one-sorted algebra $\Phi(X^0)$.

On the other hand, we can view algebra $\Phi(X \cup Y)$ as an object in the multi-sorted logic. Here, to every homomorphism $s : W(X \cup Y) \to W(X' \cup Y')$ corresponds a morphism $s_* : \Phi(X \cup Y) \to \Phi(X' \cup Y')$. For $u \in \Phi(X \cup Y)$ we have $s_* u \in \Phi(X' \cup Y')$. Let $u$ be an $X$-special formula. It is important to know for which $s$ the formula $s_* u$ is $X'$-special.

3.2. **Another characteristic of the type $Tp^H(\mu)$.** We would like to relate MT-type of a point to its LG-type.

Consider a special homomorphism $s : W(X^0) \to W(X)$ for an infinite set $X^0$ and its finite subset $X = \{x_1, \ldots, x_n\}$, such that $s(x) = x$ for each $x \in X$, i.e., $s$ is identical on the set $X$. According to the transition from $s$ to $s_*$, we have

$$s_* : \Phi(X^0) \to \Phi(X).$$

**Theorem 3.2.** *For each special homomorphism $s$, each special formula $u = u(x_1, \ldots, x_n; y_1, \ldots, y_m)$ and every point $\mu : W(X) \to H$ we have $u \in Tp^H(\mu)$ if and only if $s_* u \in LKer(\mu)$.*

*Proof.* We need one more view on a formula $u \in Tp^H(\mu)$. Given a point $\mu$, consider a set $A_\mu$ of the points $\eta : W(X^0) \to H$ defined by the rule $\eta(x_i) = \mu(x_i) = a_i$ for $x_i \in X$ and, $\eta(y)$ is an arbitrary element in $H$ for $y \in Y^0$. Denote

$$T_\mu = \bigcap_{\eta \in A_\mu} LKer(\eta).$$

It is proved [16], that a special formula $u$ belongs to the type $Tp^H(\mu)$ if and only if $u \in T_\mu$, which is equivalent to $Val_H^{X^0}(u) \supset A_\mu$.

Note that the formula $u$ of the kind

$$x_1 \equiv x_1 \wedge \ldots \wedge x_n \equiv x_n \wedge v(y_1, \ldots, y_m)$$

belongs to each $LKer(\eta)$ if the closed formula $v(y_1, \ldots, y_m)$ is satisfied in the algebra $H$. This means also that $T_\mu$ is not empty for every $\mu$.

Return to the special homomorphism $s : W(X^0) \to W(X)$ and consider the point $\mu s : W(X^0) \to H$. For $x_i \in X$ we have $\mu s(x_i) = \mu(x_i) = a_i$. Hence, the point $\mu s$ belongs to $A_\mu$.

Observe that for the formula $u = u(x_1, \ldots, x_n; y_1, \ldots, y_m)$, the formula $u(a_1, \ldots, a_n; y_1, \ldots, y_m)$ is satisfied in the algebra $H$ if the set $A_\mu$ lies in $Val_H^{X^0}(u)$. Thus, $\mu s$ belongs to $Val_H^{X^0}(u)$. By definition of $s_*$ we have that $\mu$ lies in $s_* Val_H^{X^0}(u) = Val_H^X(s_* u)$, which means that

$$s_* u \in LKer(\mu).$$

We proved the statement in one direction.

Conversely, let $s_* u \in LKer(\mu)$. Then

$$\mu \in Val_H^X(s_* u) = s_* Val_H^{X^0}(u)$$

and $\mu s \subset Val_H^{X^0}(u)$. Since the formula $u(a_1, \ldots, a_n; y_1, \ldots, y_m)$ is satisfied in $H$, then every point from the set $A_\mu$ belongs to $Val_H^{X^0}(u)$ (see also [1]). This means that the formula belongs to $Tp^H(\mu)$. $\square$

Remind that we mentioned the notion of saturated algebra. It was LG-saturation. In the Model Theory MT-saturation is defined. MT-saturation of the algebra $H$ means that for any $X$-type $T$ there is a point $\mu : W(X) \to H$ such that $T \subset Tp^H(\mu)$.

**Theorem 3.3.** *If algebra $H$ is LG-saturated then it is MT-saturated.*

*Proof.* Let algebra $H$ be LG-saturated and $T$ be $X$-MT-type correlated with $H$. We can assume that the theory $Th^{X^0}(H)$ is contained in the set of formulas $T$.

Take a special homomorphism $s : W(X^0) \to W(X)$ and pass to $s_* : \Phi(X^0) \to \Phi(X)$. Take a formula $s_* u \in \Phi(X)$ for each formula $u \in T$ and denote the set of all such $s_*u$ by $s_*T$. This set is a filter in $\Phi(X)$ containing the elementary theory $Th^X(H)$, since, if $u \in Th^{X^0}(H)$ then $s_*u \in Th^X(H)$.

Further we embed the filter $s_*T$ into the ultrafilter $T_0$ in $\Phi(X)$ which contains the theory $Th^X(H)$. By the LG-saturation of the algebra $H$ condition, $T_0 = LKer(\mu)$ for some point $\mu : W(X) \to H$. Thus, $s_*u \in LKer(\mu)$ for each formula $u \in T$. Hence, $u \in Tp^H(\mu)$ for each $u \in T$, and $T \subset Tp^H(\mu)$. This gives MT-saturation of the algebra $H$. □

We do not know whether MT-saturation implies LG-saturation. It seems that not. If it is the case, then LG-saturation of an algebra $H$ is stronger than its MT-saturation.

3.3. **Isotypeness of algebras.** The following important theorem (see 3.7 for the proof) helps to define correctly the notion of isotypeness of algebras.

**Theorem 3.4.** [19] *Let the points $\mu : W(X) \to H_1$ and $\nu : W(X) \to H_2$ be given. The equality $Tp^{H_1}(\mu) = Tp^{H_2}(\nu)$ takes place if and only if $LKer(\mu) = LKer(\nu)$ holds true.*

**Definition 3.5.** Given $X$, denote by $S^X(H)$ a set of LG-types for an algebra $H$, implemented (realized) by points in $H$. Algebras $H_1$ and $H_2$ are called isotypic if $S^X(H_1) = S^X(H_2)$ for any $X \in \Gamma^0$.

So, Theorem 3.4 implies

**Corollary 3.6.** *Algebras $H_1$ and $H_2$ in the variety $\Theta$ are isotypic if and only if they are LG-equivalent.*

According to Theorem 3.4, it doesn't matter which type is used (LG-type or MT-type).

If algebras $H_1$ and $H_2$ are isotypic then they are locally isomorphic. This means that if $A$ is a finitely generated subalgebra in $H$, then there exists a subalgebra $B$ in $H_2$ which is isomorphic to $A$ and, similarly, in the direction from $H_2$ to $H_1$.

On the other hand, local isomorphism of $H_1$ and $H_2$ does not imply their isotypeness: the groups $F_n$ and $F_m$, $m,n > 1$ are locally isomorphic, but they are isotypic only for $n = m$.

Isotypeness imply elementary equivalence of algebras, but the same example with $F_n$ and $F_m$ shows that the opposite is wrong.

Let us give some more definitions. In view of Corollary 3.6 we rephrase Definitions 2.6–3.8

**Definition 3.7.** Let $\Theta$ be a variety of algebras. We call an algebra $H \in \Theta$ **logically separable** if for any $H' \in \Theta$ isotypeness of $H$ and $H'$ implies their isomorphism.

**Definition 3.8.** The variety $\Theta$ is called **logically regular** if each free in $\Theta$ algebra $W(X)$, $X \in \Gamma^0$ is separable in $\Theta$.

**Definition 3.9.** Algebra $H \in \Theta$ is called **logically homogeneous if for every two points** $\mu : W(X) \to H$ **and** $\nu : W(X) \to H$ **the equality** $Tp^H(\mu) = Tp^H(\nu)$ **holds if and only if the points** $\mu$ **and** $\nu$ **are conjugated by an automorphism** $\sigma \in Aut(H)$, **i.e.,** $\nu = \mu\sigma$.

Logical homogenity means also, that for any point $\mu : W(X) \to H$ its $Aut(H)$ orbit is a definable set in respect to the type $LKer(\mu)$. We have also algebraic homogeneity which means that $Ker(\mu) = Ker(\nu)$ implies that the points $\mu$ and $\nu$ are $Aut(H)$ conjugated.

**Definition 3.10.** The variety $\Theta$ is called **logically perfect** if each free in $\Theta$ algebra $W(X)$, $X \in \Gamma^0$ is logically homogeneous.

The following theorem is valid:

**Theorem 3.11. If the variety $\Theta$ is logically perfect, then it is logically regular.**

*Proof.* Let the variety $\Theta$ be logically perfect and $W = W(X)$ a free in $\Theta$ algebra of the rang $n$, $X = \{x_1, \ldots, x_n\}$. Rewrite $W = H = <a_1, \ldots, a_n>$, where $a_1, \ldots, a_n$ are free generators in $H$. Let $H$ and $G \in \Theta$ be isotypic.

Take $\mu : W(X) \to H$ with $\mu(x_i) = a_i$. We have $\nu : W(X) \to G$ with $T_P^H(\mu) = T_P^G(\nu)$, $\nu(x_i) = b_i$, $B = <b_1, \ldots, b_n>$. The algebras $H$ and $B$ are isomorphic by the isomorphism $a_i \to b_i$, $i = 1, \ldots, n$.

Indeed, $T_P^H(\mu) = T_P^G(\nu)$ implies $LKer(\mu) = LKer(\nu)$ and, hence, $Ker(\mu) = Ker(\nu)$. This gives the needed isomorphism $H \to B$.

Let us prove that $B = G$. Let $B \neq G$ and there is a $b \in G$ which doesn't lie in $B$.

Take a subalgebra $B' = <b, b_1, \ldots, b_n>$ in $G$ and a collection of variables $Y = \{y, x_1, \ldots, x_n\}$ with $\nu' : W(Y) \to G$, $\nu'(y) = b$, $\nu'(x_i) = \nu(x_i) = b_i$, $i = 1, \ldots, n$.

We have $\mu' : W(Y) \to H$ with $T_P^H(\mu') = T_P^G(\nu')$. Let $\mu'(y) = a'$, $\mu'(x_i) = a'_i$, $i = 1, \ldots, n$. Let the algebras $H' = <a', a'_1, \ldots, a'_n>$ and $B' = <b, b_1, \ldots, b_n>$ be isomorphic.

Further we work with the equality $LKer(\mu') = LKer(\nu')$. Take a formula $u \in LKer(\mu)$ and pass to a formula $u' = (y \equiv y) \wedge u$. The point $(b_1, \ldots, b_n)$ satisfies the formula $u$ and, hence, the point $\nu'$ satisfies $u'$. Therefore, the point $\mu'$ satisfies $u'$ as well, and $u' \in LKer(\mu')$.

Take now a point $\mu'' : W(X) \to H$ setting $\mu''(x_i) = a'_i$, $i = 1, \ldots, n$. The point $\mu'$ satisfies the formula $u'$ if and only if the point $\mu''$ satisfies $u$. Hence, $LKer(\mu) = LKer(\mu'')$. Therefore, the point $\mu''$ is conjugated with the point $\mu$ by some isomorphism $\sigma$. Thus, the point $<a'_1, \ldots, a'_n>$ is a basis in $H$ and $a' \in <a'_1, \ldots, a'_n>$. This contradicts with $b \notin <b_1, \ldots, b_n>$. So, $B = G$ and $H$ and $G$ are isomorphic.
□

3.4. **Problems.** It seems that logical regularity of a variety $\Theta$ doesn't imply its logical perfectness, which leads to the problem

**Problem 14. Find a logically regular but not logically homogeneous variety $\Theta$. In particular, consider this problem for different varieties of groups and varieties of semigroups.**

Let us give some examples.

The variety of all groups, the variety of abelian groups and the variety of all nilpotent groups of class $n$ are logically perfect, and, hence, logically regular [5], [6], [7], [18], [19]. The variety of all semigroups and the variety of inverse semigroups are logically regular and we need to check whether they are logically perfect.

The next problem goes in parallel with the previous one:

**Problem 15. What can be said about logical regularity and logical perfectness for the variety of all solvable groups of the derived length $n$.**

Let us point some questions motivated by the example.

**Problem 16. Let $\Theta$ be a classical variety $Com - P$, a variety of commutative and associative algebras over a field $P$. The problem is to verify its logical regularity and logical perfectness.**

This is one of the problems related to logical geometry in the classical $\Theta$. The other questions on this topic will be listed later.

The pointed problem leads to the

**Problem 17. Let $S$ be a semigroup and $P$ a field, both logically homogeneous. Whether it is true that the semigroup algebra is logically homogeneous as well.**

We consider the problems of logical regularity and logical perfectness for the following varieties:

**Problem 18. The variety $Ass - P$ of associative algebras over a field $P$.**

**Problem 19. The variety $Lee - P$ of Lee algebras over a field $P$.**

**Problem 20. The variety of $n$-nilpotent associative algebras.**

**Problem 21. The variety of $n$-nilpotent Lee algebras.**

Other problems:

**Problem 22. Let $H_1$ and $H_2$ be two finitely generated isotypic algebras. Are they always isomorphic?**

**Problem 23. Let $G_1$ and $G_2$ be two finitely generated isotypic groups. Are they always isomorphic?**

**Problem 24. Let $H_1$ be a finitely generated algebra and $H_2$ is an isotypic to it algebra. Is $H_2$ also finitely generated?**

**Problem 25. Give various examples of non-commutative isotypic but not isomorphic groups. In particular, two free groups of infinite rank.**

The conditions when isotypeness does not imply isomorphism is equally interesting as the conditions when it does.

3.5. **MT-type definable sets.** We distinguish MT-definable sets and LG-definable sets. Let $A$ be a set of points in the space $Hom(W(X), H)$. This set is $X$-LG-definable if there exist a set of formulas $T$ in the algebra of formulas $\Phi(X)$ such that $A = T_H^L$. This also means that the point $\mu$ lies in $A$ if and only if it satisfies each formula $u \in T$. In other words, $T \subset LKer(\mu)$. We have also
$$T_H^L = \bigcap_{u \in T} Val_H^X(u).$$

In the case of Model Theory we take an $X$-MT-type for $T$. We set: a point $\mu \in A$ if $T \subset Tp^H(\mu)$. Let us explain this inclusion. The point $\mu$ satisfies a special formula $u = u(x_1, \ldots, x_n; y_1, \ldots, y_m)$ if the closed formula $u(a_1, \ldots, a_n; y_1, \ldots, y_m)$ holds in $H$. The inclusion $T \subset Tp^H(\mu)$ means that the point $\mu$ satisfies each formula $u \in T$.

Denote by $T_H^{L_0}$ the set of all points $\{\mu : W(X) \to H | T \subset Tp^H(\mu)\}$. We call a set $A$ MT-definable if there exists an $X$-type $T$ such that $A = T_H^{L_0}$. In the sequel we will show that each MT-definable set is LG-definable. Hence, a problem arises:

**Problem 26. Build a set $A$ which is LG-definable but is not MT-definable.**

First we need to clarify some details. Take a special morphism $s : W(X^0) \to W(X)$ identical on the set $X \subset X^0$, $X \in \Gamma^0$. We have also $s_* : \Phi(X^0) \to \Phi(X)$. Define a set of formulas $s_*T = \{s_*u | u \in T\}$.

**Theorem 3.12.** The equality $T_H^{L_0} = (s_*T)_H^L$ holds for for every $X$-type $T$.

*Proof.* Let $\mu \in T_H^{L_0}$. Then $T \subset T_P^H(\mu)$ and every formula $u \in T$ is contained in $T_P^H(\mu)$. Besides, $s_*u \in LKer(\mu)$ and $\mu \in Val_H^X(u)$. We have $\mu \in \bigcap_{u \in T} Val_H^X(u) = (s_*T)_H^L$.

Let now $\mu \in (s_*T)_H^L$. Then for every $u \in T$ we have $\mu \in Val_H^X(s_*u)$ and $s_*u \in LKer(\mu)$. Hence, $u \in Tp^H(\mu)$. This gives $T \subset T_P^H(\mu)$ and $\mu \in T_H^{L_0}$. $\square$

We see that every MT-definable set is LG-definable. The opposite statement is a problem.

Let us formulate the previous problem in another way.

**Problem 27. Whether there exists an elementary set in the LG-theory which cannot be represented as a definable set for MT-types.**

Consider now a case when algebra $H$ is logically homogeneous and $A$ is $Aut(H)$-orbit over the point $\mu : W(X) \to H$. We have $A = (LKer(\mu))_H^L$. The equality $LKer(\mu) = LKer(\nu)$ holds if and only if a point $\nu$ belongs to $A$. The same condition is needed for the equality $Tp^H(\mu) = Tp^H(\nu)$. Now, $\nu \in (Tp^H(\mu))_H^{L_0}$ by the definition of $L_0$. Thus, $A = (Tp^H(\mu))_H^{L_0}$. We proved that the orbit $A$ is MT-definable and LG-definable.

Let us present another formula for $T_H^{L_0}$. We have

$$T_H^{L_0} = \bigcap_{u \in T} Val_H^{X_0}(u).$$

Here $u = u(x_1, \ldots, x_n; y_1, \ldots, y_m)$ is a special formula in $T$ and $Val_H^{X_0}(u)$ is a set of points $\mu : W(X) \to H$ satisfying the formula $u$.

We proceed from fixed $H \in \Theta$ and $X = \{x_1, \ldots, x_n\} \in \Gamma^0$. Let us continue the definition of a Galois correspondence. Let $A$ be a subset in the points space $Hom(W(X), H)$. Let us relate to it an $X$-type $T$ by the rule

$$T = A_H^{L_0} = \bigcap_{\mu \in A} Tp^H(\mu).$$

It is checked that a special formula $u$ belongs to $T$ if and only if $A \subset Val_H^{X_0}(u)$. In other words, each point $\mu \in A$ satisfies a special formula.

This definition allows to consider Galois closures of the types $T$ and sets $A$.

Let us remind that we distinguished two full sub-categories $K_\Theta(H)$ and $LK_\Theta(H)$ in the category $Set_\Theta(H)$. Let us take one more subcategory there and denote it by $L_0K_\Theta(H)$. In each object $(X, A)$ of this category the set $A$ is an $X$-MT-type definable set. The category $L_0K_\Theta(H)$ is a full subcategory in $LK_\Theta(H)$.

**Problem 28.** [1] **Let algebras $H_1$ and $H_2$ be isotypic. Whether it is true that the categories $L_0K_\Theta(H_1)$ and $L_0K_\Theta(H_2)$ are isomorphic.**

We know that if algebras $H_1$ and $H_2$ are isotypic then the categories $LK_\Theta(H_1)$ and $LK_\Theta(H_2)$ are isomorphic. We need to check whether such isomorphism implies isomorphism of the corresponding subcategories.

We had defined Galois correspondence and, thus, we can speak about Galois closures for $A \subset Hom(W(X), H)$ and for $X$-type $T$ in MT. Namely, let us take $A^{L_0} = T$. Here $u \in T$ if $Val_H^{X_0}(\mu) \supset A$. Now,

$$A^{L_0 L_0} = T^{L_0} = \bigcap_{u \in T} Val_H^{X_0}(u).$$

Take $T_H^{L_0} = A$ for $X$-type $T$. Then $\mu : W(X) \to H$ lies in $A$ if $T \subset Tp^H(\mu)$, $u \in T^{L_0 L_0}$ if and only if $u \in A_H^{L_0}$, $Val_H^{X_0}(u) \supset A$.

It is clear that the equality $A^{L_0 L_0} = A$ means that the set $A$ is MT-definable. Analogously, the equality $A^{LL} = A$ means that the set $A$ is $LG$-definable. Thus, along with the Problem 27 we come up with the following problem

**Problem 29. Whether it is true that the equality $A^{L_0 L_0} = A$ is equivalent to $A^{LL} = A$.**

This fact seems not to be true in general. However, this is true for logically noetherian algebras $H$.

**Definition 3.13. Two algebras $H_1$ and $H_2$ are called $MT$-equivalent if $T_{H_1}^{L_0 L_0} = T_{H_2}^{L_0 L_0}$ for any $X \in \Gamma$ and $X$-type $T$.**

**Problem 30. If $H_1$ and $H_2$ are $MT$-equivalent, then the categories of definable sets $L_0K_\Theta(H_1)$ and $L_0K_\Theta(H_2)$ are isomorphic.**

We named the problems which arise naturally in the system of notions under consideration. We had not estimated the difficulty of these problems: some of them are difficult while others just need a straightforward check. We hadn't touched this issue.

---

[1] See Theorem 3.14 for a solution of Problems 26–28

In the conclusion, we compare, once again, different approaches to the notion of a definable set in the affine space $Hom(W(X), H)$. We have fixed the variety of algebra $\Theta$, an algebra $H \in \Theta$ and the finite set $X = \{x_1, \ldots, x_n\}$.

In the affine space $Hom(W(X), H)$ consider subsets $A$, whose points have the form $\mu : W(X) \to H$. Each point $\mu : W(X) \to H$ has a classical kernel $Ker(\mu)$, a logical kernel $LKer(\mu)$ and a type $(Tp^H(\mu))$. Correspondingly, we have three different geometries: algebraic geometry ($AG$), logical geometry ($LG$), and the model-theoretic geometry ($MTG$).

For $AG$ consider a system $T$ of equations $w \equiv w'$, $w, w' \in W(X)$. For $LG$ we take a set of formulas $T$ in the algebra of formulas $\Phi(X)$. For $MTG$ we proceed from an $X$-type $T$. In all these cases the set can be infinite.

Now,

• A set $A$ in $Hom(W(X), H)$ is definable in $AG$ if there exists $T$ in $W(X)$ such that $T'_H = A$, where

$$T'_H = \{\mu|\ T \subset Ker(\mu)\}$$

• A set $A$ in $Hom(W(X), H)$ is definable in $LG$ if there exists $T$ in $\Phi(X)$ such that $T^L_H = A$, where

$$T^L_H = \{\mu|\ T \subset LKer(\mu)\}$$

• A set $A$ in $Hom(W(X), H)$ is definable in $MTG$ if there exists an $X$-type $T$ such that $T^{L_0}_H = A$, where

$$T^{L_0}_H = \{\mu|\ T \subset Tp^H(\mu).\}$$

Besides that we have three closures: $T''_H$ for $AG$, $T^{LL}_H$ for $LG$, and $T^{L_0 L_0}_H$ for $MTG$. In the reverse direction the Galois correspondence for each of three cases above is as follows

$$T = A^{L_0}_H = \bigcap_{\mu \in A} Ker\mu).$$

$$T = A^{L_0}_H = \bigcap_{\mu \in A} LKer\mu).$$

$$T = A^{L_0}_H = \bigcap_{\mu \in A} Tp^H(\mu).$$

Correspondingly, we distinguish three types of equivalence relations on algebras from the variety $\Theta$.

Algebras $H_1$ and $H_2$ are algebraically equivalent if
$$T''_{H_1} = T''_{H_2}.$$
Algebras $H_1$ and $H_2$ are logically equivalent if
$$T^{LL}_{H_1} = T^{LL}_{H_2}.$$
Algebras $H_1$ and $H_2$ are $MT$-equivalent if
$$T^{L_0 L_0}_{H_1} = T^{L_0 L_0}_{H_2}.$$

3.6. **Addendum.** Return to the transition $L^0$ and check that this transition indeed determines a Galois correspondence. Let the variety $\Theta$ and $H \in \Theta$ be fixed. Let $X$ be an infinite set and $X = \{x_1, \ldots, x_n\}$ be the subset in $X$. Take the affine space $Hom(W(X), H)$ and let $\mathbb{P}$ be the system of all subsets in $Hom(W(X), H)$. Let $\mathbb{Q}$ denote the system of all $X$-types $T$ in the algebra $\Phi(X^0)$.

For $T \subset \mathbb{Q}$ we have
$$T^{L_0}_H = A = \{\mu | \ T \subset Tp^H(\mu).\}$$
Correspondingly, for $A \subset \mathbb{P}$ we have
$$T = A^{L_0}_H = \bigcap_{\mu \in A} Tp^H(\mu).$$
$T$ is an $X$-type in $\Phi(X^0)$, $T \subset \mathbb{Q}$ and consists of all $X$-special formulas such that $A \subset Val^{X_0}_H(u)$.

Check now conditions of Galois correspondence. Let $T_1 \subset T_2$, check that $T^{L_0}_{1H} \supset T^{L_0}_{2H}$. Denote $T^{L_0}_{1H} = A$ and $T^{L_0}_{1H} = B$. Let $\mu \in B$. Then $T_2 \subset Tp^H(\mu)$. Since $T_1 \subset T_2$, then $Tp^H(\mu) \supset T_1$ and $\mu \in A$. We have $B \subset A$.

Let now $A \subset B$. Check that
$$A^{L_0}_H = T_1 \supset B^{L_0}_H = T_2.$$
Let $u \in T_2$. Then $Val^{X_0}_H(u) \supset B$. Then $Val^{X_0}_H(u) \supset A$. Hence, $u \in T_1$. Thus, $T_2 \subset T_1$.

It remains to show that $A \subset A^{L_0 L_0}_H$ and $T \subset T^{L_0 L_0}_H$. For $A$ in $Hom(W(X), H)$ denote $T = A^{L_0}_H$. We have
$$A^{L_0 L_0}_H = T^{L_0} = \bigcap_{u \in T} Val^{X_0}_H(u).$$
By definition, $A$ lies in each $Val^{X_0}_H(u)$ and thus
$$A \subset A^{L_0 L_0}_H.$$
Check that for any $X$-type $T$ we have $T \subset T^{L_0 L_0}_H$. Let $A = T^{L_0}_H$. We know that $A \subset Val^{X_0}_H(u)$ for every $u \in T$. Besides that $T^{L_0 L_0}_H$ consists of all formulas $v$ such that $A \subset Val_H^{X_0}(v)$. Hence, every $u \in T$ lies in $T_H^{L_0 L_0}$ and

$$T \subset T_H^{L_0 L_0}.$$

**3.7. LG-types and MT-types.** Now, for the sake of completeness and for the aim to make picture clear and transparent we give a proof of the principal Theorem 3.17 of G.Zhitomiskii (see [19] for the original exposition). This fact is essentially used in the proof of Theorem 3.14. We hope this will help to reveal ties between two approaches to the idea of a type of a point: the one-sorted model theoretic approach and the multi-sorted logically-geometric approach. Note that the proofs are sometimes different from that of [19].

First of all, let us prove the following important fact which clarifies some of the problems (Problems 26, 27, 28) mentioned above.

**Theorem 3.14.** *Let $A \subset Hom(W(X), H)$. The set $A$ is LG-definable if and only if $A$ is MT-definable.*

*Proof.* As we know from Theorem 3.12 every $MT$-definable set is $LG$-definable. Prove the opposite.

We will use the following theorem from [19]: for every formula $u \in \Phi(X)$ there exists an $X$-special formula $\widetilde{u} \in \Phi(X^0)$ such that a point $\mu : W(X) \to H$ satisfies $\widetilde{u}$ if and only if it satisfies $u$. Let now the set $T_H^L = A$ be given. Every point $\mu$ from $A$ satisfies every formula $u \in T$. Given $T$ take $T'$ consisting of all $\widetilde{u}$ which correspond $u \in T$. The points $\mu \in A$ satisfy every formula from $T'$. This means that $T'$ is a consistent set of $X$-special formulas. Thus $T'$ is an $X$-type, such that $A \subset T_H'^{L_0}$.

Let now the point $\nu$ lies in $T_H'^{L_0}$. Then $\nu$ satisfies every formula $\widetilde{u}$. Hence it satisfies every formula $u \in T$. Thus $\nu$ lies in $T_H^L = A$. This means that

$$T_H'^{L_0} = A$$

and the theorem is proved.
□

**Corollary 3.15.** *The category $LK_\Theta(H)$ of all LG-definable sets coincides with the category $L_0 K_\Theta(H)$ of all MT-definable sets.*

Beforehand, we have proved that if the algebras $H_1$ and $H_2$ are isotypic, then the categories $LK_\Theta(H_1)$ and $LK_\Theta(H_2)$ are isomorphic. Now, the same fact is true with respect to categories $L_0 K_\Theta(H_1)$ and $L_0 K_\Theta(H_2)$.

All these provide a solution of Problems 26–28. However, we did not change the original exposition in the paper, since this insight provides the ways of the development of the topic.

**Definition 3.16.** A formula $u \in \Phi(X)$ is called correct, if there exists an $X$-special formula $\widetilde{u}$ in $\Phi(X^0)$ such that for every point $\mu : W(X) \to H$ we have $u \in LKer\mu$ if and only if $\widetilde{u} \in T_p^H(\mu)$. Denote $LG^H(\mu) = LKer\mu$.

The next theorem of G.Zhitomirskii is used in the proof of Theorem 3.14.

**Theorem 3.17.** *For every $X = \{x_1, \ldots, x_n\}$ every formula $u \in \Phi(X)$ is correct.*

*Proof.* First of all, each equality $w = w'$, $w, w' \in W(X)$ is a correct formula. This follows from $\widetilde{(w = w')} = (w = w')$.

Take two correct formulas $u$ and $v$, both from $\Phi(X)$. Show that $u \wedge v$, $u \vee v$ and $\neg u$ are also correct. We have $\widetilde{u}$ and $\widetilde{v}$. Define

$$\widetilde{u \wedge v} = \widetilde{u} \wedge \widetilde{v},$$

$$\widetilde{u \vee v} = \widetilde{u} \vee \widetilde{v},$$

$$\widetilde{\neg u} = \neg \widetilde{u}.$$

By definition, we have $u \in LKer\mu$ if and only if $\widetilde{u} \in T_p^H(\mu)$ for every point $\mu : W(X) \to H$. The same is true with respect to $v$ and $\neg u$. Let $u \vee v \in LKer\mu$ and, say, $u \in LKer\mu$. Then $\widetilde{u} \in T_p^H(\mu)$, and, hence, $\widetilde{u} \vee \widetilde{v} = \widetilde{u \vee v} \in T_p^H(\mu)$. Conversely, let $\widetilde{u \vee v} = \widetilde{u} \vee \widetilde{v} \in T_p^H(\mu)$. Suppose that $\widetilde{u} \in T_p^H(\mu)$. Then $u \in LKer\mu$, that is $u \vee v \in LKer\mu$. The similar proofs work for the correctness of the formulas $u \wedge v$ and $\neg u$. In the latter case one should use the completeness property of a type: $\neg u \in T_p^H(\mu)$ if and only if $u \notin T_p^H(\mu)$.

Our next aim is to check that if the formula $u \in \Phi(X)$ is correct, then the formula $\exists x u \in \Phi(X)$ is also correct.

Beforehand, note that it is hard to define free and bounded variables in the algebra $\Phi(X)$. This is because of the multi-sorted nature of $\Phi(X)$ and the presence in it of the formulas which include operations of the type $s_*$. So, the syntactical definition of $\exists x u \in \Phi(X)$ is a sort of problem and we will proceed from the semantical definition of this formula.

Namely, a point $\mu : W(X) \to H$ satisfies the formula $\exists x u \in \Phi(X)$ if there exits a point $\nu : W(X) \to H$ such that $u \in LKer(\nu)$ and $\mu$ coincides with $\nu$ for every variable $x' \neq x$, $x' \in X$.

Indeed, a point $\mu : W(X) \to H$ satisfies $\exists x u \in \Phi(X)$ if $\mu \in Val_H^X(\exists x u) = \exists x(Val_H^X(u))$ (see Subsection 1.4). Denote the set $Val_H^X(u)$ in $Hal_\Theta^X(H) = Bool(W(X), H)$ by $A$. Then $\mu$ belongs to $\exists x A$. Using the definition of existential quantifiers in $Hal_\Theta^X(H)$ (Subsection 1.3) and the fact that $u \in LKer(\nu)$ if and only if $\nu \in Val_H^X(u)$, we arrive to the definition above.

Since $u$ is correct, there exists an $X$-special formula $\widetilde{u} \in \Phi(X^0)$,

$$\widetilde{u} = \widetilde{u}(x_1, \ldots, x_n, y_1, \ldots, y_m), \ x_i \in X, \ y_i \in Y^0 = (X^0 \setminus X),$$

such that $\widetilde{u} \in T_p^H(\mu)$ if and only if $u \in LKer(\mu)$, where $\mu : W(X) \to H$.

Define

$$\widetilde{\exists x u} = \exists x \widetilde{u}.$$

The formula $\exists x \widetilde{u}$ is not $X$-special since $x$ is bound (we assume that $x$ coincides with one of $x_i$, say $x_n$). Take a variable $y \in X^0$, such that $y$ is different from each $x_i \in X$ and $y_j \in Y^0$.

Define $\exists y \widetilde{u}_y$ to be a formula which coincides with $\exists x \widetilde{u}$ modulo replacement of $x$ by $y$. So, $\exists y \widetilde{u}_y$ has one less free variable and one more bound variable than $\exists x \widetilde{u}$.

Consider endomorphism $s$ of $W(X^0)$ taking $s(x)$ to $y$ and leaving all other variables from $X^0$ unchanged. Let $s_*$ be the corresponding automorphism of the one-sorted Halmos algebra $\Phi(X^0)$. Then $s_*(\exists x \widetilde{u}) = \exists s_*(x) s_*(\widetilde{u}) = \exists y \widetilde{u}_y$.

Define

$$\widetilde{\exists x u} = \exists y \widetilde{u}_y.$$

Thus, in order to check that $\exists x u$ is correct, we need to verify that for every $\mu : W(X) \to H$ the formula $\exists x u$ lies in $LKer(\mu)$ if and only if $\exists y \widetilde{u}_y \in T_p^H(\mu)$.

Let $\exists x u$ lies in $LKer(\mu)$. Thus, there exits a point $\nu : W(X) \to H$ such that $u \in LKer(\nu)$ and $\mu$ coincides with $\nu$ for every variable $x' \neq x$, $x' \in X$. Consider $X_y = \{x_1, \ldots, x_{n-1}, y\}$.

We have points $\mu : W(X) \to H$, $\mu' : X_y \to H$ where $\mu'(x_i) = \mu(x_i) = a_i$, and $\mu'(y)$ is an arbitrary element $b$ in $H$. We have also $\nu : W(X) \to H$ and $\nu' : X_y \to H$, where $\nu'(x_i) = \nu(x_i)$, and $\nu'(y) = \nu(x_n)$. So, $\nu$ and $\nu'$ have the same images. Denote it $(a_1, a_2, \ldots, a_{n-1}, a_n)$, $a_i \in H$, i.e., $\nu'(y) = a_n$.

Take

$$\widetilde{u}_y = \widetilde{u}(x_1, \ldots, x_{n-1}, y, y_1, \ldots, y_m),$$

Since the formula $\exists y \widetilde{u}(a_1, \ldots, a_{n-1}, b, y_1, \ldots, y_m)$ is closed for any $b$, then either it is satisfied on any point $\mu'$, or no one of $\mu'$ satisfies this formula. We can take $b = a_n$, that is $\mu' = \nu'$. Since $\nu$ and $\nu'$ have the same images, and $u$ is correct, the point $\nu'$ satisfies $\widetilde{u}_y$. Then $\nu'$ satisfies $\exists y \widetilde{u}_y$. Hence $\exists y \widetilde{u}(x_1, \ldots, x_{n-1}, y, y_1, \ldots, y_m)$ is satisfied on $\mu'$ for any $b$. This means that $\exists y \widetilde{u}_y \in T_p^H(\mu')$ for every $\mu'$. We can take $\mu'$ to be $\mu$. Then $\widetilde{\exists x u} \in T_p^H(\mu)$.

Conversely, let $\widetilde{\exists x u} \in T_p^H(\mu)$. Take a point $\nu : W(X) \to H$ such that $\nu(x_i) = \mu(x_i)$, $i = 1, \ldots, n-1$, $\nu(x_n) = \nu(y)$. We have $\widetilde{u} \in T_p^H(\nu)$. Since $\widetilde{u}$ is correct, then $u$ in $LKer(\nu)$. The points $\mu$ and $\nu$ coincide on all $x_i$, $i \neq n$. Thus $\exists u$ belongs to $LKer(\mu)$.

It remains to check that the operation $s_*$ respects correctness of formulas. Let $X = \{x_1, \ldots, x_n\}$, $Y = \{y_1, \ldots, y_m\}$, and a morphism $s : W(Y) \to W(X)$ be given. Take the corresponding $s_* : \Phi(Y) \to \Phi(X)$. Given $v \in \Phi(Y)$ consider $u = s_* v$ in $\Phi(X)$. We shall show that if $v$ is $Y$-correct then $u$ is $X$-correct.

We have $u \in LKer(\mu)$, $\mu : W(X) \to H$ if and only if $v \in LKer(\nu)$, $\nu : W(Y) \to H$ for $\mu s = \nu$. Indeed, $u = s_* v \in LKer(\mu)$ means that $\mu \in Val_H^X(s_* v) = s_* Val_H^Y(v)$ and thus, $\mu s \in Val_H^Y(v)$. Hence, for $\nu = \mu s$ we have $v \in LKer(\nu)$. Conversely, let $v \in LKer(\nu)$ and $\mu s = \nu \in Val_H^Y(v)$. We have $\mu \in s_* Val_H^Y(v) = Val_H^X(s_* v) = Val_H^X(u)$ and $u \in LKer(\mu)$.

Note that morphism $s_* : \Phi(Y) \to \Phi(X)$ is a homomorphism of boolean algebras. Suppose that $v \in \Phi(Y)$ is correct. We have

$$\widetilde{v} = \widetilde{v}(y_1, \ldots, y_m, z_1, \ldots, z_t),$$

where all $z_i$ are bound and belong to $Z = \{z_1, \ldots, z_t\}$. All free variables in $\widetilde{v}$ belong to $Y$ (it is assumed that not necessarily all variables from $Y$ occurs in $\widetilde{v}$). In this sense $\widetilde{v}$ is $Y$-special.

We will define also the formula $\widetilde{u}$ and show that in our situation $\widetilde{u} \in Tp^H(\mu)$ if and only if $\widetilde{v} \in Tp^H(\nu)$.

Consider $Z' = \{z_1', \ldots, z_t'\}$, where all $z_i'$ do not belong to $X$. Take the free algebras $W(X \cup Z')$ and $W(Y \cup Z)$. Define homomorphism $s' : W(Y \cup Z) \to W(X \cup Z')$ extending $s : W(Y) \to W(X)$ by $s'(z_i) = z_i'$. The commutative diagram of homomorphisms takes place:

$$\begin{array}{ccc} W(Y \cup Z) & \xrightarrow{s'} & W(X \cup Z') \\ {\scriptstyle s^1}\downarrow & & \downarrow{\scriptstyle s^2} \\ W(Y) & \xrightarrow{s} & W(X). \end{array}$$

Here $s^1$ and $s^2$ are special homomorphisms which act identically on $Y$ and $X$, respectively. The corresponding commutative diagram of morphisms of algebras of formulas is as follows:

$$\begin{array}{ccc} \Phi(Y \cup Z) & \xrightarrow{s'_*} & \Phi(X \cup Z') \\ {\scriptstyle s^1_*}\downarrow & & \downarrow{\scriptstyle s^2_*} \\ \Phi(Y) & \xrightarrow{s_*} & \Phi(X). \end{array}$$

This diagram is commutative due to the fact that the product of morphisms of algebras of formulas corresponds to the product of homomorphisms of free algebras. Apply the diagram to $Y$-special formula $\widetilde{v}$ which belongs to the algebra $\Phi(Y \cup Z)$. Then, $s^2_* s'_* \widetilde{v} = s_* s^1_* \widetilde{v}$. Assume that $\widetilde{u} = s'_* \widetilde{v}$. Here, $\widetilde{u}$ is an $X$-special formula, contained in the algebra $\Phi(X \cup Z')$. We need to prove that for any point $\mu : W(X) \to H$ the inclusion $\widetilde{u} \in Tp^H(\mu)$ holds if and only if $u \in LKer(\mu)$.

We use the criterion from Section 3 (Theorem 3.2): $\widetilde{u} \in Tp^H(\mu)$ if and only if $s^2_* \widetilde{u} \in LKer(\mu)$. Let us prove the latter inclusion. The similar criterion is valid for the formula $\widetilde{v}$. Since the formula $v$ is correct, then $\widetilde{v} \in Tp^H(\nu)$, where $\nu = \mu s$. Hence, $s^1_* \widetilde{v} \in LKer(\nu)$, which means that the point $\nu$ belongs to the set $Val^Y_H(s^1_* \widetilde{v})$. Since $\nu = \mu s$, then $\mu \in Val^X_H(s_* s^1_* \widetilde{v}) = Val^X_H(s^2_* s'_* \widetilde{v}) = Val^X_H(s^2_* \widetilde{u})$. This leads to the inclusion $s^2_* \widetilde{u} \in LKer(\mu)$, which gives $\widetilde{u} \in Tp^H(\mu)$.

The same reasoning in the opposite direction shows that the inclusion $\widetilde{u} \in Tp^H(\mu)$ is equivalent to that of $\widetilde{v} \in Tp^H(\nu)$.

It is worth to recall that we started from the fact $u \in LKer(\mu)$ if and only if $v \in LKer(\nu)$. But, $v \in LKer(\nu)$ because of the correctness of the formula $v$. Thus, $u \in LKer(\mu)$. Hence, the transition from $u$ to $\widetilde{u}$ guarantees the correctness of the formula $u$.

Hence, the set of all correct $X$-formulas, for various $X$, respects all operations of the multi-sorted algebra $\widetilde{\Phi}$. Since $\widetilde{\Phi}$ is generated by equalities, which are correct, the subalgebra of all correct formulas in $\widetilde{\Phi}$ coincides with $\widetilde{\Phi}$. Thus every $u \in \widetilde{\Phi}(X)$ for every $X$, is correct. □

**Theorem 3.18.** *Let the points $\mu : W(X) \to H_1$ and $\nu : W(X) \to H_2$ be given. Then*

$$Tp^{H_1}(\mu) = Tp^{H_2}(\nu)$$

*if and only if*

$$LKer(\mu) = LKer(\nu).$$

*Proof.* Let the points $\mu : W(X) \to H_1$ and $\mu : W(X) \to H_2$ be given and let $Tp^{H_1}(\mu) = Tp^{H_2}(\mu)$. Take $u \in LKer(\mu)$. Then $\widetilde{u} \in Tp^{H_1}(\mu)$ and, thus, $\widetilde{u} \in Tp^{H_2}(\nu)$. Hence, $u \in LKer(\nu)$. The same is true in the opposite direction.

Let, conversely, $LKer(\mu) = LKer(\nu)$. Take an arbitrary $X$-special formula $u$ in $Tp^{H_1}(\mu)$. Take a special homomorphism from $s : W(X^0) \to W(X)$. It corresponds the morphism $s_* : \Phi(X^0) \to \Phi(X)$. Then, using Theorem 3.2, the formula $u \in Tp^H(\mu)$ if and only if $s_*u \in LKer(\mu)$. Then $s_*u \in LKer(\nu)$. Then $u \in Tp^H(\nu)$. $\square$

Consider a simple example. Take $Y = \{y_1, y_2\}$ and $X = \{x_1, x_2, x_3\}$ and let $s$ be a homomorphism $s : W(Y) \to W(X)$. Take also variables $z$ and $z'$ and extend $s$ to $s' : W(Y \cup z) \to W(X \cup z')$ assuming $s'(z) = z'$. We have also morphism $s'_* : \Phi(Y \cup z) \to \Phi(X \cup z')$. Take an equality $w(y_1, y_2, z) \equiv w'(y_1, y_2, z)$ in $\Phi(Y \cup Z)$. Consider $\exists z(w \equiv w')$ and apply $s'_*$. We have
$$s'_*(\exists z(w \equiv w')) = \exists z'(s'w \equiv s'w').$$
Here $s'w = s'(w(y_1, y_2, z) = w(w_1, w_2, z')$, where $w_i = s(y_1) = w_i(x_1, x_2, x_3)$ and

$$s'_*(\exists z(w \equiv w') = \exists z'(w(w_1(x_1, x_2, x_3), w_2(x_1, x_2, x_3), z')$$

$$\equiv w'(w_1(x_1, x_2, x_3), w_2(x_1, x_2, x_3), z')).$$

In the conclusion one more problem which is connected with the previously named problems on isotypeness and isomorphism of free algebras.

**Problem 31.** Let two isotypic finitely-generated free algebras $H_1$ and $H_2$ and two points $\mu : W(X) \to H_1$ and $\nu : W(X) \to H_2$ be given. Let $LKer(\mu) = LKer(\nu)$. Is it true that there exists an isomorphism $\sigma : H_1 \to H_2$ such that $\mu\sigma = \nu$?

Boris Plotkin: Institute of Mathematics, Hebrew University, 91904, Jerusalem, ISRAEL

*E-mail address*: `plotkin` *at* `macs.biu.ac.il`